\documentclass[12pt]{amsart}

\usepackage{color}
\usepackage[all]{xy}
\usepackage{amssymb,amsmath}

\usepackage{mathrsfs}
\usepackage{eucal}

\usepackage{setspace}

\date{March 3, 2014}

\calclayout
\newtheorem{dummy}{anything}[section]
\newtheorem{theorem}[dummy]{Theorem}

\theoremstyle{definition}

  \newtheorem{remark}[dummy]{Remark}

\newcommand{\cS}{\mathscr S}

\newcommand{\bZ}{\mathbb Z}

\newcommand{\bQ}{\mathbb Q}

\newcommand{\bbH}{\mathbb H}

\newcommand{\bbQ}{\mathbb Q}
\newcommand{\bbR}{\mathbb R}
\newcommand{\bbS}{\mathbb S}
\newcommand{\bbZ}{\mathbb Z}
\DeclareMathOperator{\Image}{Image}

 \DeclareMathOperator{\rank}{rank}

\newcommand{\cy}[1]{\bZ/{#1}}
\newcommand{\la}{\langle}
\newcommand{\ra}{\rangle}
\newcommand{\vv}{\, | \,}
\newcommand{\bd}{\partial}

\newcommand\units[1]{(#1)^{\times}}

\newcommand\wX{\widetilde X}
\newcommand\wH{\widehat H}
\newcommand\wK{\widetilde K}
\newcommand\Zpi{\bZ\pi}
\newcommand\Qhat{\widehat\bbQ}

\newcommand\Qpi{\bbQ\pi}
\newcommand\Zhat{\widehat\bbZ}
\DeclareMathOperator\vcd{vcd}

\begin{document}

\title{Topological Spherical Space Forms}
\author{Ian Hambleton}

\address{Department of Mathematics, McMaster University,
Hamilton, Ontario L8S 4K1, Canada}

\email{hambleton@mcmaster.ca }

\thanks{Research partially supported by NSERC Discovery Grant A4000.}

\begin{abstract} 
Free  actions of finite groups on spheres give rise to  topological  spherical space forms. The existence and classification problems for space forms have a long history in the  geometry and topology of manifolds. In this article, we present a survey of some of the main results and a guide to the literature.
\end{abstract}



\maketitle

\nocite{*}

\section{Introduction}
A Clifford-Klein space form is a closed Riemannian manifold with constant sectional curvature. A classical result in differential geometry states that space forms are  just the quotients of the standard models  $\bbS^n$, $\bbR^n$ and $\bbH^n$ for spherical, Euclidean, or hyperbolic geometry by discrete groups of isometries acting freely and properly discontinuously (see \cite[Chap.~8]{Carmo:1992}). The spherical space forms have now been completely classified, answering the original question of Killing \cite{Killing:1891} from 1890.  Extensive information for both compact and non-compact space forms is contained in Wolf \cite{Wolf:2011}. 

In 1926, Hopf \cite{Hopf:1926} determined the fundamental groups of three-dimensional spherical space forms, and asked  whether every closed, oriented $3$-manifold with finite fundamental group admits a Riemannian metric of constant positive curvature. In particular, this topological question is a generalization of the famous Poincar\'e conjecture to the non-simply connected case. By 2005, both Hopf's question and the Poincar\'e conjecture were finally answered by Perelman \cite{Lott:2007}, building on the Ricci flow techniques of Hamilton \cite{Hamilton:1982}, in the process of veriflying the Thurston Geometrization conjecture (see 
\cite{Cao:2006,Kleiner:2008,Morgan:1984,Morgan:2007}).

The \emph{topological} spherical space form problem \cite{Wall:1978} generalizes Hopf's question to higher dimensions: which finite groups can act freely on spheres via homeomorphisms or diffeomorphisms ? How can one classify free (smooth or topological) finite group actions on spheres ? These questions have led to an enormous amount of research in algebraic and geometric topology over the last hundred years, and substantial progress has been made. However, interesting aspects of these questions are still unresolved, so we may expect future work in this fascinating area of mathematics.

The purpose of this survey is to give an overview of the current knowledge about  the topological spherical space form problem, and to provide a reading guide to the literature. We apologize for any errors or omissions, and note that the previous surveys listed in the references contain much more information than can be included here (see, for example: \cite{Adem:2002,Davis:1985,Hambleton:2002,Madsen:1980,Madsen:1981,Milgram:1982,Milnor:1957,Smith:1960,Wall:1978}). Some new directions are briefly described in the last section.

\section{First steps: pre-1960}\label{sec:firststeps}
The development of algebraic topology in the 1930's led to the first necessary conditions for a finite group $\pi$ to act freely on a sphere. Since $\cy 2$ is the only possibility for a free action on $S^{2k}$ (by the Lefschetz fixed point theorem), we assume from now on that our space forms are \emph{odd dimensional}. Foundational work on group actions was done by P.~A.~Smith \cite{Smith:1938,Smith:1960} in the 1940's, who showed that, for each prime $p$,  subgroups of order $p^2$ in such a group $\pi$ must be cyclic. More precisely, the elementary abelian $p$-group 
$\cy p\times \cy p$ can not act freely and simplicially on any finite simplicial complex which is homotopy equivalent to a sphere. 

We say that a finite group $\pi$ such that \emph{every subgroup of order $p^2$ is cyclic} satisfies the \emph{$p^2$-condition}. Equivalently, the $p^2$-condition means that the $p$-Sylow subgroup is cyclic, or possibly generalized quaternion (if $p=2$). It now becomes a problem in group theory to classify the finite groups $\pi$ satisfying the $p^2$-conditions for all primes.  This problem was settled   by Zassenhaus \cite{Zassenhaus:1935} and Suzuki \cite{Suzuki:1955}: the quotient of $\pi$ by the maximal normal odd order subgroup $O(\pi)$ belongs to one of six types:
\begin{enumerate}
\item[I.] cyclic
\item[II.] generalized quaternion $2$-group
\item[III.] binary tetrahedral
\item[IV.] binary octahedral
\item[V.] $SL_2(p)$, $p\geq 5$ prime,
\item[VI.] $TL_2(p)$, $p\geq 5$ prime.
\end{enumerate}
The groups of type I-IV are solvable, and the groups of type V-VI are non-solvable. The notation $TL_2(p)$ indicates a certain extension of $SL_2(p)$ (see \cite[p.~103]{Thomas:1971} for more detail). Zassenhaus showed that if $\pi$ admits a free orthogonal action on a sphere, then every subgroup of order $pq$, $p$ and $q$ primes, must be cyclic. Conversely, if a solvable group $\pi$ satisfies all the $pq$-conditions, then it admits a free orthogonal representation. 

After the introduction of group cohomology by  Hurewicz, Artin and Tate, an alternative characterization was given: the finite groups satisfying the $p^2$-conditions for all primes are precisely the finite groups with periodic cohomology \cite[Chap.~XII, \S 11]{Cartan:1999}. In this formulation,  cup product with an  element $g\in \wH^n(\pi;\bZ)$   gives an isomorphism $\wH^k(\pi;A) \xrightarrow{\approx} \wH^{n+k}(\pi;A)$, $k \in \bZ$, for any $\pi$-module $A$. If $n>0$, such elements $g$ are in bijection with \emph{generators} of the cyclic group $\wH^n(\pi;\bZ)\cong \cy{|\pi|}$. The degree $n$ is called a cohomological \emph{period}, and it turns out that every period is a multiple of an even minimal period, denoted $2d(\pi)$.

Thomas and Wall \cite{Thomas:1971} defined a \emph{$(\pi,n)$-polarized space}, for $n\geq 3$,  to consist of a pointed CW-complex $X$, dominated by a finite complex, together with an isomorphism $\pi_1(X, x_0) \cong \pi$ and a homotopy equivalence of the universal covering $\wX$ to the sphere $S^{n-1}$. Here are two key observations:

\begin{enumerate}
\addtolength{\itemsep}{0.2\baselineskip}
\item the (polarized) homotopy type of a $(\pi,n)$-polarized space is determined, via the first $k$-invariant, by a generator of $\wH^n(\pi;\bZ)$, and
\item each $(\pi,n)$-polarized space $X$ determines an obstruction
 (see Wall  \cite{Wall:1967})
$$\theta(X) = \theta(g) \in \wK_0(\Zpi)$$
 such that $X$ is homotopy equivalent to a finite complex if and only if $\theta(g) = 0$. 
\end{enumerate}
A smooth or topological free $\pi$ action on $S^{n-1}$ leads to a closed manifold quotient $X = S^{n-1}/\pi$.  In this case, $X$ is a $(\pi,n)$-polarized space, whose $k$-invariant  $g\in \wH^n(\pi;\bZ)$  has the property $\theta(g) = 0$. Finding such a space is the first basic obstruction to the existence problem. We will discuss some of the delicate number theory issues involved in Section \ref{sec:finiteness}. 

In the late 1950's, Milnor's dramatic discoveries revolutionized algebraic and geometric topology. For the space form problem, the discovery of exotic spheres \cite{Milnor:1956} and his update of Hopf's 1926 paper about $3$-manifolds with finite fundamental group had the most striking impact. 
\begin{theorem}[Milnor \cite{Milnor:1957}]\label{thm:2p}
 If $\pi$ is a finite group which acts freely by homeomorphisms on a manifold which is a mod 2 homology sphere, then every subgroup of order $2p$ in $\pi$ is cyclic.
\end{theorem}
After this the topological space form problem suddenly had two versions, depending on whether one considers smooth actions on (possibly exotic) homotopy spheres, or topological actions on the standard sphere. Moreover, Milnor had identified a new group theoretic obstruction to existence beyond the $p^2$-conditions obtained by P.~A.~Smith. For a free $\pi$-action on  $S^{n-1}$ by homeomorphisms, $\pi$ must satisfy the \emph{$2p$-conditions} for all primes. Equivalently, $\pi$ cannot contain any dihedral subgroups or every element of order two must be central. Milnor pointed out that the first unsolved case after this new restriction is the non-abelian group of order 21.

 For actions on the $3$-sphere, $\pi$ must be $4$-periodic.
In the same paper, Milnor  \cite{Milnor:1957}
set the stage for future work by listing the finite groups which satisfy all the $p^2$ and $2p$ conditions, and have period four. Among the groups in Milnor's list, which did not appear in Hopf's list (arising from orthogonal spherical space forms $\bbS^3/\Gamma$, with $\Gamma \subset SO(4)$), were an interesting family:
$$
Q(8m, k, l) = \{x, y, z\vv x^2 = (xy)^2 = y^{2m}, z^{kl} = 1, xzx^{-1} = z^r, yzy^{-1} = z^{-1}\}
$$
where $8m$, $k$, $l$ are pairwise relatively prime positive integers,  with $r \equiv -1 \pmod k$, and $r \equiv +1 \pmod l$.  If $l=1$ we write $Q(8m,k)$, and if $m$ is odd we label these group \emph{type IIM}. This family of $4$-periodic groups  has played a pivotal role in all the subsequent developments. 

\section{Swan complexes: 1960-1970}
The next big step was taken in 1960 by Swan, who established a converse to the necessary conditions from P.~A.~Smith theory for simplicial actions.
\begin{theorem}[Swan \cite{Swan:1960}] Let $\pi$ be a finite group with periodic cohomology. For $n \geq 4$, the set of (polarized) homotopy types of  $(\pi,n)$-polarized spaces is in bijection with the set of generators  $g \in \wH^n(\pi;\bZ)$. 
\end{theorem}
From now on, we will usually refer to $(\pi,n)$-polarized spaces as \emph{Swan complexes}. Notice that $\pi$ is only assumed to satisfy the $p^2$-conditions, so dihedral subgroups present no obstruction to the existence of Swan complexes.
\begin{proof}
Here is an outline of the proof of this fundamental result. Swan first showed that every periodicity generator $g \in \wH^n(\pi;\bZ)$ for the cohomology of $\pi$ determines a periodic projective resolution
$$P_*(g) : 0 \to \bZ \to P_{n-1} \to P_{n-2} \to \dots \to P_1 \to P_0 \to \bZ \to 0$$
meaning an exact sequence of $\Zpi$-modules, where each $P_i$ is a finitely-generated projective module over the group ring $\Zpi$ and $\bZ$ denotes the integers with trivial $\pi$-action. Given $g$, the resolution $P_*(g)$ is unique up to a chain homotopy equivalence which is the identity on $\bZ$ at both ends.
The geometric realization of $P_*(g)$ is carried out by showing that the resolution $P_*(g)$ may be assumed to agree in degrees $i \leq 2$ with the cellular chain complex $C_*(\wK)$ of a finite $2$-complex $K$ with fundamental group $\pi$. Then an inductive cell-attaching argument (due to Milnor) produces a finite-dimensional CW-complex with a free cellular action of $\pi$. The quotient of this CW-complex by the action of $\pi$ is the desired Swan complex.
\end{proof}

In order to produce a finite Swan complex, it is necessary to analyze the effect of changing the generator. 
 If $g' =r\cdot g\in \wH^n(\pi;\bZ)$ is another generator, then the comparison theorem for resolutions \cite[III.6.1]{Mac-Lane:1995} gives a chain map $\phi\colon P_*(g) \to P_*(g')$
$$\xymatrix{
0 \ar[r] &\bZ \ar[r] \ar[d]^r& P_{n-1}\ar[d]^\phi\ar[r] & P_{n-2}\ar[d]^\phi \ar[r] & \dots \ar[r] & P_1 \ar[d]^\phi\ar[r] & P_0\ar[r]\ar[d]^\phi &\bZ\ar[r]\ar@{=}[d] & 0\cr
0 \ar[r] &\bZ \ar[r] & P'_{n-1} \ar[r] & P'_{n-2} \ar[r] & \dots \ar[r] & P'_1 \ar[r] & P'_0\ar[r] &\bZ\ar[r] & 0
}$$
inducing multiplication by $r$  on the left-hand module $\bZ$.  Note that $r$ is relatively prime to the order of $\pi$, since $g$ and $g'$ are both generators.

For each generator, Swan now defined the finiteness obstruction $\theta(g) \in \wK_0(\Zpi)$ by the Euler characteristic formula
$$ \theta(g) = \sum_{i = 1}^{n-1} (-1)^i[P_i] \in \wK_0(\Zpi)$$
 where $[P_i]$ denotes the class of  $P_i$ in  the projective class group.
 He showed that $\theta(g) = 0$ if and only if $g$ can be represented by a periodic resolution of finitely-generated \emph{free} $\Zpi$-modules. 
 Moreover, if $g'=r\cdot g$ is another generator, then $\theta(g') = \theta(g) + [\la N, r\ra]$, where 
 $\la N, r\ra$ denotes the ideal in $\Zpi$ generated by the sum $N$ of all the group elements of $\pi$, and an integer $r$ relatively prime to  $|\pi|$. 

 
 These  projective ideals are usually called \emph{Swan modules}, and the subgroup of $\wK_0(\Zpi)$ generated by all the Swan modules (as $r$ varies) is called the
 \emph{Swan subgroup} $T(\pi)$ of the projective class group.   It follows that a periodic group $\pi$ admits a finite $(\pi, n)$-complex if and only if the image 
 $$ \theta_n(\pi):= [\theta(g)] = 0 \in \wK_0(\Zpi)/T(\pi).$$ 
 In other words, 
   the indeterminacy in the finiteness obstruction given by changing the generator is precisely described by $T(\pi)$.

For the minimal period dimension
  $n = 2 d(\pi)$, we obtain an invariant $\theta(\pi):=[\theta(g)]$ depending only on $\pi$. This turned out to be very difficult to calculate.
  On the other hand, if we are willing to increase the dimension by splicing $k$ copies of a given periodic projective resolution $P_*(g)$, we obtain a resolution $P_*(g^k)$ of length $nk-1$. This algebraic operation corresponds to the geometric join of $k$ copies of $\wX(g)$. From the definition, it follows that $\theta(g^k) = k\cdot \theta(g) \in \wK_0(\Zpi)$. In another fundamental paper, Swan \cite{Swan:1960a} proved that $\wK_0(\Zpi)$ is a finite abelian group. 
Hence  the finiteness obstruction can be eliminated by choosing $k$ equal to the exponent of the projective class group. In summary:
\begin{theorem}[Swan \cite{Swan:1960}] Let $\pi$ be a finite group with periodic cohomology. Then, for some integer $n = 2kd(\pi)$, there exists a finite $(n-1)$-dimensional Swan complex $X(\pi)$.
\end{theorem}
 This result settles the ``stable existence" problem for Swan complexes, but existence in the period dimension $n= 2d(\pi)$ turned out to require difficult calculations in number theory.

 \section{Semicharacteristic obstructions and $3$-manifolds}
 
 In the late 1960's,
Wall asked whether Milnor's result on the necessity of the $2p$-conditions could be interpreted
 using the new theory of non-simply connected  surgery.
Ronnie Lee \cite{Lee:1973a} answered  this question in 1973 by
defining a ``semicharacteristic"  obstruction for the
problem (see also \cite{Davis:1983,Pardon:1980}).
As well as recovering the previous result of
Milnor, 
the semicharacteristic shows that none of the groups
$Q(8m,k,l)$, $m$ even, defined at the end of  Section \ref{sec:firststeps}  can act freely 
and topologically on a sphere $S^{8s+3}$, for any $s \geq 0$.
Later Thomas \cite{Thomas:1988} observed that this  
also eliminates the family of groups
$O(48, 3^{k-1},l)$ in Milnor's list from  dimensions $8s + 3$ because groups of this type always
 contain a subgroup isomorphic to
$Q(16,3^{k-1},1)$. 

In particular,  the difference between Hopf's list of orthogonal space form groups and Milnor's list of possible finite fundamental groups of closed $3$-manifolds could now be resolved by eliminating a single family of $4$-periodic groups, namely  the type IIM groups $Q(8m,k,l)$, with $m$ odd.  However, each of these groups   contains a subgroup of the form
$Q(8p,q)$, for distinct odd primes $p$ and $q$, so it  is enough to eliminate  the groups $Q(8p,q)$. 
Despite claims in the literature (see \cite{Thomas:1978,Thomas:1979}), only Perelman's results 
more than 30 years later gave the first proof that $Q(8p,q)$ is not a $3$-manifold group.

 \section{Poincar\'e complexes}
 
In 1967, Wall \cite{Wall:1967} defined Poincar\'e complexes as CW complexes with the homology properties of closed manifolds, and it was clear that every $(\pi,n)$-polarized space $X$ was a Poincar\'e complex in this sense.  
Surgery theory was developed by Browder \cite{Browder:1972}, Kervaire-Milnor \cite{Kervaire:1963},  Novikov, Sullivan and Wall \cite{Wall:1999} to study the existence and uniqueness of manifold structures on Poincar\'e duality spaces. For surgery on topological manifolds one also needs Kirby-Siebenmann \cite{Kirby:1977}.

 The \emph{structure set}, $ \cS(X)$, is the set of equivalence classes of pairs $(M,f)$, where $M^m$ is a closed, topological or smooth $m$-manifold, and $f\colon M \to X$ is (simple) homotopy equivalence. The existence problem is to show $\cS(X) \neq \emptyset$, and surgery theory sets up a two step process for the solution. In the first step, one must find a (topological) vector bundle $\xi \searrow X$ which is fibre homotopy equivalent to the Spivak normal fibre space  $\nu_X$ for $X$ (see  \cite{Spivak:1967}, \cite[Chap.~3]{Wall:1967}). This bundle  is a candidate for the stable normal bundle of a manifold homotopy equivalent to $X$. It is classified by a lift of the classifying map $\nu_X\colon X \to BG$ to $BTOP$ or $BO$. If such a bundle exists, we say that $\nu_X$ is \emph{reducible}, and one may vary the choice of reduction by an element in $[X, G/TOP]$ or $[X, G/O]$. A reduction for $\nu_X$ provides a \emph{degree one normal map} or \emph{normal invariant}
 $$\xymatrix{\nu_M \ar[r]^b\ar[d] & \xi\ar[d]\cr
 M \ar[r]^f & X
 }$$
 which is well-defined up to \emph{normal cobordism}. The second step in the surgery process is to evaluate an obstruction
 $$\sigma(f, b) \in L_m(\Zpi)$$
of the normal invariant in an algebraically defined \emph{surgery obstruction group}. For $m \geq 5$, this obstruction vanishes if and only if $(f, b)$ is normally cobordant to a (simple) homotopy equivalence. In other words, the outcome of the two step process (if successful) is the existence of a closed TOP or smooth manifold homotopy equivalent to $X$. 
 
If $\cS(X) \neq \emptyset$, $m\geq 5$,  and we fix a smooth manifold structure on $X$ as a base-point, the surgery exact sequence (see Wall \cite[Chap.~10]{Wall:1999})
 $$  [\Sigma(X), G/O] \xrightarrow{\sigma_{m+1}(X)} L_{m+1}(\Zpi) \longrightarrow \cS(X) \longrightarrow [X, G/O] \xrightarrow{\sigma_m(X)} L_m(\Zpi)$$
  relates the classification of manifolds which are simple homotopy equivalent to $X$, or the determining the set $\cS(X)$,  to the calculation of the surgery obstruction maps $\sigma_{m+1}(X)$ and $\sigma_m(X)$. Here 
  $\Sigma(X) = (X\times I)/(X \times \bd I)$, and \emph{exactness} is interpreted in the context of pointed sets.  In particular, $\sigma_m(X) \colon [X, G/O] \to L_m(\Zpi)$ is \emph{not} a homomorphism in  general, although both domain and range are abelian groups. There is a similar sequence  for topological manifolds, with $[X, G/TOP]$ instead of $[X, G/O]$ and the obvious change in the definition of the structure set. 
  
   For a suitable $H$-space structure on $G/TOP$ (see \cite{Kirby:1977}, \cite{Madsen:1979}, \cite{Ranicki:1992}), the surgery exact sequence becomes a sequence of abelian groups and homomorphisms. In particular,  the structure set $\cS_{TOP}(X)$ inherits a natural abelian group structure in the topological case (see Ranicki  \cite{Ranicki:1992}).
 
 In 1971 Ronnie Lee \cite{Lee:1971} and Petrie \cite{Petrie:1971} independently showed how surgery theory could be used to prove a new existence result.
 
 \begin{theorem}[Lee, Petrie] There exist free smooth actions of odd order non-abelian metacyclic groups on spheres.
\end{theorem}
In particular, these groups $\pi=\cy p \rtimes \cy q$ are solvable but do not satisfy the $pq$-conditions, so the actions constructed here gave the first examples of  \emph{non-linear} free actions on spheres.
 
 \section{Surgery theory and stable existence}
 In the 1970's, Madsen, Thomas and Wall combined ideas from group cohomology and representation theory, with new methods in algebraic topology and hermitian forms, to completely solve the stable existence problem for topological space forms.
 \begin{theorem}[\cite{Madsen:1976, Madsen:1983c}] Let $\pi$ be a finite group   satisfying the $p^2$ and $2p$-conditions, for every prime $p$. Then $\pi$ acts freely and smoothly on $S^{n-1}$, for any $n\equiv 0 \pmod{2e(\pi)}$.
\end{theorem}
In the statement, $e(\pi)$ is the Artin-Lam induction exponent \cite{Curtis:1990}. This is the least positive integer $e$ such that $e\cdot \chi$ is expressible as a sum of rational representation induced from cyclic subgroups of $\pi$, for each $\chi \in R_\bbQ(\pi)$. For the periodic groups $e(\pi)= d(\pi)$ or $e(\pi) = 2d(\pi)$, depending on $\pi$, and  the dimensional result is best possible,  except possibly for groups containing a type IIM subgroup $Q(8m,k)$.


We can only outline the ingredients in the proof. The key idea is to choose the homotopy type and normal invariant of the Swan complex $X(\pi)$ so that its coverings $X(\rho)$, for proper subgroups $\rho\subset \pi$ agree with orthogonal spherical space forms (whenever this makes sense). To achieve this goal,  one can use the fact that both group cohomology $\wH^*(\pi;\bZ)$ and the normal invariants $[X(\pi), G/O]$ are described by (generalized) cohomology theories, and hence computable by induction and restriction in terms of the Sylow subgroups.  For the normal invariant, this is major new ingredient based on the infinite loop space structure of $G/O$ (see Madsen and Milgram \cite{Madsen:1979} for references and more information). For the surgery step, Wall \cite{Wall:1976}, \cite[Chap.~13B]{Wall:1999} had shown that 
\begin{enumerate}
\item  the surgery obstruction groups $L_*(\Zpi)$, for $\pi$ finite, are computable in the sense of Dress induction theory  by restriction to the $2$-hyperelementary and the odd $p$-elementary subgroups of $\pi$, and
\item the change in the surgery obstruction by varying the normal invariant is detected by restriction to the $2$-Sylow subgroups of $\pi$.
\end{enumerate}
For periodic groups satisfying the $2p$-conditions, all the $2$-hyperelementary and odd $p$-elementary groups admit free orthogonal representations. The main existence result of Madsen, Thomas and Wall now follows directly from these general facts and the careful preparation of the surgery problem. Here is a summary technical statement (obtained by combining \cite[\S\S 2-3]{Madsen:1976} and \cite[Theorem 2]{Madsen:1983c}).

\begin{theorem}[Madsen-Thomas-Wall]
Let $\pi$ be a finite group with periodic cohomology. There exist finite $(\pi,n)$-polarized Swan complexes $X = X(\pi)$ such that the covering spaces $X(\rho)$ are homotopy equivalent to closed manifolds, for each $\rho\subseteq\pi$ which has a fixed-point free orthogonal representation. In addition, $X$ has a smooth normal invariant which restricts to the normal invariant of an orthogonal spherical space form for the $2$-Sylow covering $X(\pi_2)$.
\end{theorem}
In this formulation, $\pi$ is not assumed to satisfy the $2p$-conditions and this extra generality is useful for other applications. However, as mentioned above, if $\pi$ has no dihedral subgroups then the surgery obstruction automatically vanishes for a suitable smooth normal invariant, implying that $\pi$ acts freely and smoothly on a sphere.

\section{The finiteness obstruction}\label{sec:finiteness}
The remaining part of the topological spherical space form problem, at least for existence, concerns actions in the period dimension $n = 2d(\pi)$. Here the first major difficulty is to produce a \emph{finite} $(\pi,2 d(\pi))$-polarized Swan complex.  The problem is made more challenging by the fact that the projective class group $\wK_0(\Zpi)$ containing the finiteness obstructions requires (in general) all the hyperelementary subgroups for a computation using Dress induction. An effective method of calculation was described by Wall \cite{Wall:1979} in 1979, based on a Mayer-Vietoris sequence arising from the arithmetic square:
$$\xymatrix{\Zpi\ar[r]\ar[d]&\Qpi\ar[d]\cr
\Zhat\pi \ar[r] & \Qhat\pi
}$$
In this sequence
$$ \dots\to K_1(\Zhat\pi) \oplus K_1(\Qpi) \to K_1(\Qhat\pi) \xrightarrow{\bd}
\wK_0(\Zpi) \to \wK_0(\Zhat\pi) \oplus \wK_0(\Qpi)\to  \dots$$
the finiteness obstruction $\theta(g)$ for a periodic projective resolution 
$$P_*(g) : 0 \to \bZ \to P_{n-1} \to P_{n-2} \to \dots \to P_1 \to P_0 \to \bZ \to 0$$
lies in the image of the boundary map $\bd\colon K_1(\Qhat\pi) \to \wK_0(\Zpi)$. More precisely, we have the formula
$$ \bd(\widehat\Delta(g)) = \theta(g) \in \wK_0(\Zpi)$$
where $\widehat\Delta(g) \in K_1(\Qhat\pi)$ is the \emph{idelic Reidemeister torsion} of the resolution, or its associated $(\pi, n)$-polarized Swan complex. 
The definition of  the idelic Reidemeister torsion in  \cite[\S 9]{Wall:1979} follows the prescription in Milnor \cite{Milnor:1966} for the torsion of a chain complex of free modules with based homology. In this case, $P_*(g) \otimes R$ becomes a chain complex of free $R$-modules, for $R = \bbQ$ or $ \Zhat$, since  projective modules over $\Zpi$ are locally free, by Swan \cite{Swan:1960}. The complex $P_*(g)$ is the pull-back of based chain complexes over $\Zhat\pi$ and $\Qhat\pi$ via glueing automorphisms over $\Qhat\pi$, leading to $\widehat\Delta(g) \in K_1(\Qhat\pi)$.

The challenge was to actually carry out the calculations in specific cases. As already explained, the main examples to consider are the $4$-periodic type IIM groups. The breakthrough here was made in 1981 by 
Milgram, who first showed that the finiteness obstruction  $\theta(g)$ lies in a more tractable subgroup $D(\Zpi) \subseteq \wK_0(\Zpi)$ of the projective class group (defined in \cite[Chap.~6, p.~234]{Curtis:1987}), and then obtained a surprising result:
\begin{theorem}[{Milgram \cite{Milgram:1981}}]   Among the periodic groups $\pi = Q(8p,q)$, for distinct odd primes $p$, $q$,  there exist groups with non-zero finiteness obstruction,
and groups which can act freely and topologically on $S^{8s+3}$, for $s >0$.
\end{theorem}
After this,  further concrete calculations were made by Milgram \cite{Milgram:1985}, \cite{Davis:1985}, and Madsen \cite{Madsen:1980}. The methods were re-organized and extended by  Bentzen and Madsen \cite{Bentzen:1983, Bentzen:1984, Bentzen:1987}. At present, it is still not known exactly which groups $Q(8p,q)$ have non-zero finiteness obstruction. On the other hand, assuming the Generalized Riemann Hypothesis, there are infinitely prime pairs $(p,q)$ for which $Q(8p,q)$ acts freely on $S^{8s+3}$, $s >0$ (see Bentzen \cite[6.6]{Bentzen:1987}). For the range $pq < 2000$, Bentzen lists 
$$(p,q) = (3,313), \, (3, 433), \, (3, 601)\ \text{ and } \ (17,103)$$
as the only pairs for which $Q(8p,q)$ admits a free action on $S^{8s+3}$. These calculations were made possible by using reduced norms to identify the idelic torsion in $K_1(\Qhat\pi)$ with units in the $l$-adic completions of the centre fields of $\Qpi$ in the Wedderburn decomposition. 

The situation for more general type IIM groups $Q(8p,q) \times \cy r$ is even more difficult, because the finiteness obstruction is not detected on the subgroups. For example, for $\pi=Q(24,313) \times \cy 7$ the finiteness obstruction is non-zero, but $\theta(Q(24, 313) ) = 0$.

\section{Space forms in the period dimension}

In 1983,  Madsen \cite{Madsen:1983a} translated the existence problem for a free $Q(8p,q)$-action  on a sphere $S^{8s+3}$, $s>0$, into explicit number theoretic conditions. In this sense, we have  a complete solution to the topological spherical space form problem for these groups. The work cited above on the finiteness obstruction shows that effective calculations are possible in principle, but they involve detailed knowledge of units in algebraic number fields. 

In order to give the flavour of the number theory involved, we need some definitions (see \cite[p.~195]{Madsen:1983a}). For any integer $r$, let $\eta_r = \zeta_r + \zeta_r^{-1}$, where $\zeta_r = e^{2\pi i/r}$ is a primitive root of unity. Let $A = \bZ[\eta_p, \eta_q]$ and let $F = \bQ(\eta_p, \eta_q)$ denote its fraction field. There are reduction maps
$$\Phi_A \colon F^{(2)} \to \units{A/pq}= \units{A/p} \oplus \units{A/q}$$
and
$$\varphi_A\colon F^{(2)} \to \units{A/4A}$$
where $F^{(2)}$ denotes the elements in $F$ with even valuation at all finite primes in $A$.
\begin{theorem}[Madsen {\cite[Theorem B]{Madsen:1983a}}]\label{thm:madsen} 
The group $Q(8p,q)$ acts freely on $S^{8s+3}$, for $s > 0$, if and only if the following conditions
are satisfied:
\begin{enumerate}
\item $(1, -1) \in \Image(\Phi_A)$;
\item $(\eta_p-2, \eta_q-2) \in \Image(\Phi_A \vv \ker \varphi_A)$;
\item $(2,2) \in \Image (\Phi_A)$.
\end{enumerate}
\end{theorem}
The first and third conditions imply that the ``top" component of the finiteness obstruction is zero.
The mysterious elements $(\eta_p-2, \eta_q-2) \in  \units{A/p} \oplus \units{A/q}$ arise from the idelic Reidemeister torsion of a suitable Swan complex, and the reduction map gives the ``top" component of the surgery obstruction. A crucial technical device used to identify explicit elements in algebraic $K$ and $L$-theory in terms of units are \emph{Galois module homomorphisms}, as introduced by Fr\"ohlich \cite{Frohlich:1983}. The proof of this existence result uses the full calculation of surgery obstruction groups \cite{Wall:1976}, \cite{Hambleton:2000}.   Madsen's statement  included the assumption that the $k$-invariant was ``almost linear", but this turned out to be a necessary condition (see  \cite{Hambleton:1986a}).  
\begin{remark} The vanishing of the high-dimensional obstructions for $\pi = Q(8p,q)$ also  implies that $\pi$ acts freely on some integral homology $3$-sphere. Conversely, the non-vanishing of the high-dimensional obstructions shows that $\pi$ is not a $3$-manifold fundamental group. 
\end{remark}

\section{Topological Euclidean space forms}
A free $\pi$-action on a sphere $S^{n-1}$ produces by ``coning" a topological $\pi$-action on Euclidean space $\bbR^{n}$. One can ask: ``which finite groups can act semi-freely on Euclidean space $\bbR^n$ with one isolated fixed point ?". More generally, ``which finite groups can act semi-freely on some Euclidean space $\bbR^{n+k}$ with (non-empty) fixed set a standardly embedded $\bbR^k$ ?" By one-point compactification, we obtain the corresponding questions for spheres, and then generalize to allow more isotropy:
\begin{enumerate}
\item[A.] Which finite groups admit a semi-free topological action on $S^{n+k}$,  with  fixed set a standard $S^k$, $k \geq 0$ ?
\item[B.] Which finite groups $\pi$ can act topologically on $S^{n+k}$, freely on the complement of a $\pi$-invariant standard sphere $S^k$, $k \geq 0$ ?
\item[C.] Which infinite discrete groups $\Gamma$ can act freely, properly discontinuously and co-compactly on
$S^m \times \bbR^n$ ? 
\end{enumerate}
The complement $S^{n+k} - S^{k} \approx S^{n-1} \times \bbR^{k+1}$ is homotopy equivalent to a sphere. Therefore, any finite group $\pi$ acting as in problem (A) or  (B) , and any finite subgroup of $\Gamma$ in problem (C),  must be a periodic group (satisfying the $p^2$-conditions for all primes $p$). 
Problem (C) is a natural extension of (B) to a new world of infinite discrete group actions.  Examples arise from the geometry of pseudo-Riemannian space forms (see Kulkarni \cite{Kulkarni:1978,Kulkarni:1981}). 
Wall \cite[p.~518]{Wall:1979} noted that in this setting the Farrell-Tate cohomology \cite{Farrell:1977} would necessarily be periodic (e.g.~if $\vcd\Gamma < \infty$), and asked whether $\Gamma$ could contain any (finite) dihedral subgroups.

\smallskip
Problem (A)  has been completely solved (reduced to number theory). Here is a sample result, for semi-free actions on $(\bbR^n, 0)$ or $(S^n, S^0)$.
\begin{theorem}[Hambleton-Madsen {\cite[Theorem C]{Hambleton:1986}}] The group $Q(8p,q)$ acts semi-freely on $(\bbR^{4l},0)$,  for $l \geq 0$, if and only if $(\eta_p-2, \eta_q-2) \in \Image(\Phi_A \vv \ker \varphi_A)$ modulo squares.
\end{theorem}
In the statement the ranges of both reduction maps used in Theorem \ref{thm:madsen} are taken modulo squares of units. This number theoretic condition is much easier to satisfy, and Bentzen \cite{Bentzen:1987} made a thorough comparison with the spherical case. For example, there are infinitely many pairs $(p,q)$ which satisfy the condition, and many examples where the finiteness obstruction of $Q(8p,q)$ is non-zero. In other words, these examples have no invariant free $\pi$-action on a sphere $S^{4l-1}$  in the complement of the isolated fixed point $\{0\}$. The extension of surgery theory for finite groups to dimension four by Freedman \cite{Freedman:1984}, \cite[\S 11.10, p.~231]{Freedman:1990} allows the statements in \cite{Hambleton:1986} to be extended to actions on $(\bbR^4, 0)$. One-point compactification produces interesting orientation-preserving examples of semi-free topological actions on $S^4$ fixing two points, by finite groups which are not subgroups of $SO(5)$. Another general  result on problem (A) is that every $\pi$ with such a semi-free action on $(S^{n+k}, S^k)$ also acts semi-freely on $(S^n, S^0)$. Moreover, every such semi-free action ``desuspends" to an action on $(S^{n+1}, S^1)$, but not necessarily to $(S^n, S^0)$ (see \cite[p.~804]{Hambleton:1986}). The desuspension result was also proved in \cite{Anderson:1983}.

Problems (B) and (C) have led to further research, but they are still very far from solution. In \cite{Hambleton:1991} these problems were studied using new methods in bounded topology, and the role of the $2p$-conditions was understood. The relevant surgery obstruction groups were computed in \cite{Hambleton:1993}. It was shown that dihedral groups can occur for (B), but only in the presence of a specific non-trivial action on the invariant sub-sphere $S^k\subset S^{n+k}$.

\begin{theorem}[Hambleton-Pedersen {\cite[7.11]{Hambleton:1991}}] 
Let $V$ be a linear representation of the dihedral group $\pi = D_{2p}$, for $p$ an odd prime. There is a topological $\pi$-action on a sphere $S^{n+k}$, which is free on the complement of a standard $\pi$-invariant sub-sphere, and given by the unit sphere $S(V)$ on the sub-sphere, if and only if the representation $V$ has at least two $\bbR_{-}$ summands.
\end{theorem}
The $\bbR_{-}$ representation is the sign representation of $D_{2p}$ via its projection onto $\cy 2$.
A similar statement \cite[8.3]{Hambleton:1991} was proved in the setting of problem (C), showing exactly when the infinite discrete groups $\Gamma = \bbZ^r \rtimes_\alpha D_{2p}$ can act freely, properly and co-compactly on some product $S^m \times \bbR^n$, expressed in terms of the twisting map $\alpha\colon D_{2p} \to GL_r(\bbZ)$.

The classification of $\pi$-actions on $(\bbR^{n+k}, \bbR^k)$ or $(S^{n+k}, S^k)$ up to equivariant homeomorphism is very challenging, even for the  linear actions provided by direct sums $V\oplus W$ of a free $\pi$-representation $V$ with an arbitrary $\pi$-representation $W$. This is the famous \emph{non-linear similarity} problem of de Rham (see \cite{Hambleton:1995,Hambleton:2005}). A lot of  work remains to be done on this problem for non-cyclic periodic groups $\pi$.
We remark that the version of problem (C) without assuming co-compactness has been solved (see 
\cite{Adem:2001,Connolly:1989}). 

\section{Further topics}
The spherical space form problem concerns free actions of finite groups on spheres. Here are two recent directions related to these periodic or \emph{rank one} groups. Recall that the \emph{rank} of a finite group $\pi$ is the largest integer $k$ such that $\pi$ contains an elementary abelian $p$-group $(\cy p)^k$ for some prime $p$.
\begin{enumerate}
\item[A.] Which finite groups can act freely on a product $S^n \times S^n$ of spheres of the same dimension ?
\item[B.] Which finite groups can act on a sphere $S^n$ with periodic isotropy subgroups ?
\end{enumerate}
The first problem has group theoretic obstructions, but the complete answer is far from understood. For example,  if $\pi$ acts freely on any product of spheres, then $\rank\pi \leq 2$. If $\pi$ acts freely on $S^n\times S^n$, for some $n$, then $\pi$ cannot contain the alternating group $A_4$ as a subgroup (see Oliver \cite{Oliver:1979}). This implies that no finite simple group can act freely on $S^n \times S^n$.  In the late 1970's, Stein \cite{Stein:1979} proved that $\pi = D_{2p} \times D_{2p}$ can act freely and smoothly on $S^{4k+3} \times S^{4k+3}$, for $k \geq 0$, even though the dihedral group can not act freely on a single sphere. 
\begin{theorem}[Hambleton \cite{Hambleton:2006}]
Let $\pi = \pi_1 \times \pi_2$ be a product of finite periodic groups. Then $\pi$  acts freely and smoothly on some product
$S^n \times S^n$.
\end{theorem}
Further work on problem (A) is contained in \cite{Hambleton:2009,Hambleton:2010}, including existence results for the extra-special $p$-groups of order $p^3$, for $p$ an odd prime. Free actions of finite groups on arbitrary products of spheres have been extensively studied (see Adem and Smith \cite{Adem:2001}).

If a finite group acts on a sphere with periodic isotropy, then $\rank\pi \leq 2$. Orthogonal representation spheres $S(V)$ provide a rich collection of examples, but there are many rank two finite groups whch do not admit any linear action with periodic isotropy. The first unknown case was $\pi = S_5$, the symmetric group on 5 letters. Here is a recent analogue of Swan's theorem for this case.
\begin{theorem}[Hambleton-Pamuk-Yal{\c{c}}{\i}n \cite{Hambleton:2013}] Let $\pi = S_5$. Then there exists a finite $\pi$-CW complex $X$ with periodic isotropy such that $X \simeq S^n$.
\end{theorem}
Much work remains to be done on problem (B), to identify the group theoretic and topological obstructions (see \cite{Hambleton:2011,Kulkarni:1982,Kwasik:1988}).  Recall that, by P.~A.~Smith theory, the
 fixed sets of $p$-power subgroups are $\cy p$-homology spheres.
The work of Jones \cite{Jones:1971} and Oliver \cite{Oliver:1975} points toward the richness of the classification problem for finite group actions on spheres.
 In addition, Milnor's $2p$-conditions on normalizer quotients of isotropy subgroups also appear in this setting as obstructions to smooth or topological actions. In future work, it will be very interesting to understand how the structure of the fixed sets interacts with the family of isotropy subgroups, in order to compare linear and non-linear actions on spheres.

\bigskip
The following list of references is organized approximately by relevance to the sections headings above. 
\providecommand{\bysame}{\leavevmode\hbox to3em{\hrulefill}\thinspace}
\providecommand{\MR}{\relax\ifhmode\unskip\space\fi MR }
\providecommand{\MRhref}[2]{%
  \href{http://www.ams.org/mathscinet-getitem?mr=#1}{#2}
}
\providecommand{\href}[2]{#2}
\renewcommand\refname{Surveys}

\end{document}